\documentclass[reqno, 12pt]{amsart}
\usepackage[letterpaper,hmargin=1in,vmargin=1in]{geometry}
\usepackage{amsmath, amssymb, amsthm, verbatim, url}
\usepackage{hyperref} 

\newcommand \N {\mathbb{N}}
\newcommand \R {\mathbb{R}}

\newcommand \Oh {\mathcal{O}}

\newcommand \eps {\varepsilon}

\newcommand \Def {\stackrel{\textrm{def}}=}

\DeclareMathOperator \supp {supp}

\DeclareMathOperator \WFh {WF_{\textit{h}}}
\DeclareMathOperator \Op {Op}
\DeclareMathOperator \Id {Id}

\newtheorem*{thm}{Theorem}

\theoremstyle{definition}

\numberwithin{equation}{section}

\title{Extending cutoff resolvent estimates via propagation of singularities}
\author{Kiril Datchev}
\email{datchev@math.mit.edu}
\address{Mathematics Department, Massachusetts Institute of Technology, Cambridge, MA 02139}

\begin{document}

\begin{abstract}
We use a gluing method developed in joint work with Andr\'as Vasy to show that polynomially bounded cutoff resolvent estimates at the real axis imply the same estimates, up to a constant factor, in a neighborhood of the real axis.
\end{abstract}

\maketitle

\section{Introduction}

Let 
\[
P = -h^2\Delta + V(x), \qquad R(\lambda) = (P-\lambda)^{-1},
\]
where $V \in C_0^\infty(\R^n)$. Suppose $\supp V \subset \{|x| < R_0\}$, and take $E>0$ and $\chi \in C_0^\infty(\R^n)$ with $\chi(x) =1$ near $\{|x| \le R_0 + 5\}$. We show that polynomial semiclassical estimates for the cutoff resolvent $\chi R(E+i0)\chi$ imply the same estimates, up to a constant factor, for the meromorphic continuation of $\chi R(\lambda)\chi$ to $\lambda$ near $E$. Cutoff functions (or at least weights) are needed to define the limit $R(E+i0) = \lim_{\eps \to 0^+} R(E+i\eps)$ and the continuation to $\lambda$ near $E$ because $E \in [0,\infty)$, which is the essential spectrum of $P$.
\begin{thm}
Suppose there exist $h_0>0$ and $a\colon (0,h_0] \to (0,\infty)$  satisfying $1 \le a(h)\le h^{-N}$ for some $N \in \N$, such that
\begin{equation}\label{e:resinit}
\|\chi R(E+i0) \chi\|_{L^2(\R^n) \to L^2(\R^n)} \le a(h), \qquad 0<h\le h_0.
\end{equation}
Then there exist $C, h_1>0$ such that the meromorphic continuation of $\chi R(\lambda) \chi$ from $\{\arg \lambda > 0\}$ to $\{\arg \lambda <0\}$ obeys
\begin{equation}\label{e:resfin}
\|\chi R(\lambda) \chi\|_{L^2(\R^n) \to L^2(\R^n)} \le C a(h), \qquad  |\lambda - E| \le \frac 1{Ca(h)}, \,\,0<h\le h_1.
\end{equation}
\end{thm}
The function $a(h)$ for which the estimate \eqref{e:resinit} is satisfied is related to dynamical properties of the trapped set of the Hamiltonian $p = |\xi|^2 + V(x)$:
\[
K_E = \{\rho \in p^{-1}(E); \,\exists C_\rho>0, \,\forall t \in \R, \,\left|\exp(tH_p)\rho\right| \le C_\rho\}.
\]
Here $H_p = 2 \xi \cdot \nabla_x - \nabla V\cdot \nabla_\xi$ is the Hamiltonian vector field associated to $p$. We have \eqref{e:resinit} with $a(h) = C/h$ for all $E$ in a neighborhood of $E_0$ if and only if $K_{E_0} = \varnothing$, and in this case \eqref{e:resfin} is well known and no better bound holds (but see \cite[Proposition 3.1]{nsz} for an estimate which covers larger range of $\lambda$).

The most interesting case of the Theorem is when $a(h) = C\log(1/h)/h$. In \cite{bbr}, Bony-Burq-Ramond show that if $K_E \ne \varnothing$, then $\sup_{E'}\|\chi R(E'+i0) \chi\| \ge \log(1/h)/(Ch)$, where the supremum is of $E'$ in a neighborhood of $E$ (although in general the optimal $a(h)$ is $\exp(C/h)$; see Burq \cite{Burq:Lower}, and in this general case the Theorem does not apply). The estimate \eqref{e:resinit} with $a(h) = C\log(1/h)/h$ is true in many situations when the trapping is hyperbolic: see  Burq \cite{bur:smoothing}, Christianson \cite{chr:noncon, c08}, Nonnenmacher-Zworski \cite{nz1, nz2} and Wunsch-Zworski \cite{wz}. In particular it is true if $K_E$ consists of a single hyperbolic orbit. In the case of a single orbit which is degenerately hyperbolic, Christianson-Wunsch \cite{cw} give examples where \eqref{e:resinit} holds for $a(h) = h^{-k}$ where $k>1$ depends on the degeneracy.

Similar results are proved using different methods by Burq \cite{b98} and Vodev \cite{v,vodev2} (who obtain a continuation with estimates in an exponentially small region, but without any assumptions on $K_E$) and Christianson \cite{c08, chr:mrl} (who proves the result when $a(h) = C\log(1/h)/h$ in the case where $K_E$ is a single hyperbolic orbit). Those papers also give applications of such resolvent estimates to wave decay (see also \cite{lebeau} for a similar result in the setting of the damped wave equation).

In the present paper we use the gluing method via propagation of singularities developed in collaboration with Vasy in \cite{d-v1}; this is perhaps the simplest application of that method. The same proof holds with only minor modifications when $V$ is replaced by a metric perturbation or an obstacle or when $\chi$ is noncompactly supported but suitably decaying. As in \cite{d-v1}, it can similarly treat suitably decaying perturbations or asymptotically hyperbolic manifolds in the sense of \cite{v1,v2}.  In particular, in the case $a(h) = C \log(1/h)/h$ above, a rescaling of parameters in the Theorem implies the conclusion of \cite[Corollary 2.3]{c08} in a more general trapping situation, namely we replace the assumption that there is one trapped orbit with \eqref{e:resinit}. This should have as a consequence the wave decay asserted in \cite[Theorem 2]{c08} under the same more general hypotheses, with only minor modifications to the proof in that paper. In the interest of simplicity we do not pursue these generalizations and applications here.

I am very grateful to Vesselin Petkov and Nicolas Burq for suggesting this problem as an application of the gluing method of \cite{d-v1}, and for helpful discussions and comments. Thanks also to the anonymous referee for several suggestions and corrections. Finally, I am grateful for the support of a National Science Foundation postdoctoral fellowship, and for the hospitality of the Universit\'e Paris Nord, Villetaneuse.

\section{Proof of Theorem}

If the cutoff $\chi$ were not present, \eqref{e:resfin} would follow directly from the resolvent identity
\begin{equation}\label{e:resid}
R(\lambda) = R(E)(\Id + (E-\lambda)R(E))^{-1}.
\end{equation}
In order to use this identity, we introduce a complex absorbing barrier function which suppresses the effects of infinity and removes the need for cutoffs in the estimates. Indeed, take $W \in C^\infty(\R^n;[0,\infty))$ with $W(x) = 0$ near $|x| \le R_0 + 4$ and with $W(x)=1$ near $|x| \ge R_0 + 5$. Put 
\[
P_W = P - iW, \qquad R_W(\lambda) = (P_W - \lambda)^{-1}, 
\]
and note that the essential spectrum of $P_W$ is $-i+[0,\infty)$ (since $P_W$ is a relatively compact perturbation of $-h^2\Delta - i$), so that $R_W(\lambda)$ is meromorphic (and, as we will see, actually holomorphic) for $\lambda$ near $E$, without needing to be multiplied by cutoffs. We will show that \eqref{e:resinit} implies
\begin{equation}\label{e:resinitw}
\|R_W(E)\| \le C a(h),
\end{equation}
and that
\begin{equation}\label{e:resfinw}
\| R_W(\lambda)\| \le C a(h), \qquad  |\lambda - E| \le \frac 1{Ca(h)},
\end{equation}
implies \eqref{e:resfin}.  The implication \eqref{e:resinitw} $\Longrightarrow$ \eqref{e:resfinw} follows from  \eqref{e:resid}. Here and below all operator norms are $L^2(\R^n) \to L^2(\R^n)$, and all function norms are $L^2(\R^n)$. The large constant $C$ may change from line to line, and $h \in (0,h_{\textrm{max}}]$ where $h_{\textrm{max}}$ may change from line to line.

We prove \eqref{e:resinit} $\Longrightarrow$ \eqref{e:resinitw} and \eqref{e:resfinw} $\Longrightarrow$ \eqref{e:resfin} using two reference operators, where the perturbation (in this case the potential $V$) is removed. Put 
\[
P_0 = -h^2\Delta, \qquad P_{W,0} = P_0 - iW,
\]
and define $R_{0}(\lambda) = (P_0 - \lambda)^{-1}$ and $R_{W,0}(\lambda) = (P_{W,0}-\lambda)^{-1}$. These resolvents obey the usual nontrapping bounds (see e.g.  \cite[Proposition 3.1]{nsz}):
\begin{equation}\label{e:infstan}
\|\chi R_{0} (\lambda) \chi\|  + \|R_{W,0} (\lambda)\|\le \frac C h, \qquad  |\lambda - E| \le \frac h C. 
\end{equation}

\begin{proof}[Proof that \eqref{e:resinit} $\Longrightarrow$ \eqref{e:resinitw}]
Take $\chi_K \in C^\infty(\R)$ with $\chi_K(r) = 1$ near $r \le R_0 + 2$ and $\chi_K(r) = 0$ near $r \ge R_0 + 3$, and put $\chi_\infty = 1 - \chi_K$. Let
\[
F = \chi_K(|x| - 1) R(E + i0) \chi_K(|x|) + \chi_\infty(|x|+1)R_{W,0}(E) \chi_\infty(|x|).
\]
 Now put
\[\begin{split}
(P_W - E)F 
&= \Id + [P,\chi_K(|x| - 1)] R(E+i0) \chi_K(|x|) + [P,\chi_\infty(|x|+1)]R_{W,0}(E)\chi_\infty(|x|) \\
&\Def  \Id + A_K + A_\infty.
\end{split}\]
These errors are not small (we only have $\|A_K\| \le   C a(h)h$ and $\|\chi A_\infty \chi\| \le C$), but solving them away using $F$ we obtain errors which we can control using propagation of singularities. In fact, using $A_K^2 = A_\infty^2 = 0$ we obtain
\[(P_W-E)F(\Id - A_K - A_\infty + A_KA_\infty) = \Id -A_\infty A_K + A_\infty A_K A_\infty.\]
Using $\chi_K(|x|)A_K = 0$, \eqref{e:resinit} and \eqref{e:infstan} we find that
\[\|\chi F(\Id - A_K - A_\infty + A_KA_\infty)\chi\| \le C a(h),\]
so the conclusion follows from
\begin{equation}\label{e:ohinf}
\| A_\infty A_K\| = \Oh(h^\infty).
\end{equation}
But this is a consequence of propagation of singularities along nontrapping bicharacteristics. Indeed, recall that the semiclassical wavefront set of a function $u \in L^2(\R^n)$ with $\|u\| \le h^{-N}$, denoted $\WFh u$, is defined as follows: a point $\rho \in T^*\R^n$ is not in $\WFh u$ if there exists $a \in C^\infty(T^*\R^n)$, bounded together with all derivatives and with $a(\rho) \ne 0$, such that $\|\Op(a)u\| = \Oh(h^\infty)$.  Define the free backward bicharacteristic at $\rho =(x,\xi) \in T^*\R^n$ by 
\[\gamma_\rho^- = \{(x + 2t\xi,\xi); t \le 0\}.\]
We will use propagation of singularities in the following form: if $f \in L^2(\R^n)$ has compact support and $\|f\| \le  h^{-N}$ for some $N \in \N$, then for all $\rho \in T^*\R^n$, 
\begin{equation}\label{e:propsingrww}
\rho \in \WFh \left(R_{W,0}(E)f\right) \Longrightarrow \gamma_\rho^- \cap \WFh f \ne \varnothing.
\end{equation}
See for example \cite[Lemma 5.1]{d-v1} for a proof. If further $\gamma_\rho^- \subset T^*\{|x|>R_0\}$ (so that $\gamma_\rho^- \cap T^*\supp V = \varnothing$ and hence $\gamma_\rho^-$ is also the backward bicharacteristic for the Hamiltonian $p = |\xi|^2 + V$), then similarly
\begin{equation}\label{e:propsingr}
\rho \in \WFh \left(R(E+i0) f\right) \Longrightarrow \gamma_\rho^- \cap \WFh f \ne \varnothing.
\end{equation}
See for example \cite[Lemma 2]{d} for a proof. We now use this to show that $\WFh A_\infty A_K f$ is empty for any $f \in L^2(\R^n)$ with $\|f\| = 1$, from which \eqref{e:ohinf} follows.  To see this, write
\[
A_\infty A_K f = [P,\chi_\infty(|x|+1)]R_{W,0}(E)  [P,\chi_K(|x| - 1)] R(E+i0) \chi_K(|x|) f.
\]
Then any $\rho = (x,\xi) \in \WFh A_\infty A_K f$ has 
\begin{equation}\label{e:prop1}
\rho \in T^* \supp d \chi_\infty(|\cdot|+1) \subset T^*\{R_0 + 1 < |x| < R_0 + 2\},
\end{equation}
and by \eqref{e:propsingrww} we know that $\gamma_\rho^-$ must contain a point 
\begin{equation}\label{e:prop2}
\rho' \in T^* \supp d \chi_K(|\cdot| - 1) \subset T^*\{R_0 + 3< |x| < R_0 + 4\}.
\end{equation}
Since $\rho' = (x+2t\xi,\xi)$ for some $t < 0$, and since \eqref{e:prop1} and \eqref{e:prop2} imply that $|x+2t\xi| > |x|$, it follows that $(x + 2t\xi) \cdot \xi<0$ and in particular $\gamma_{\rho'}^- \subset T^*\{|x|>R_0+3\}$. Because of this we may apply \eqref{e:propsingr} to conclude that $\gamma_{\rho'}^-$ must contain a point 
\[
\rho'' \in T^* \supp \chi_K(|\cdot|) \subset T^*\{|x| <  R_0 + 3\},
\]
which is impossible because $\gamma_{\rho'}^- \subset \{|x|>R_0+3\}$. This shows that $\WFh A_\infty A_K f = \varnothing$, from which \eqref{e:ohinf} follows.\end{proof}

\begin{proof}[Proof that \eqref{e:resfinw} $\Longrightarrow$ \eqref{e:resfin}]
We use the same $\chi_K$ and $\chi_\infty$ as in the previous proof, but we redefine $F$, $A_K$, and $A_\infty$ as follows:
\[
F = \chi_K(|x| - 1) R_W(\lambda) \chi_K(|x|) + \chi_\infty(|x|+1)R_{0}(\lambda) \chi_\infty(|x|).
\]
Now put
\[\begin{split}
(P - \lambda)F  &= \Id + [P,\chi_K(|x| - 1)] R_W(\lambda) \chi_K(|x|) + [P,\chi_\infty(|x|+1)]R_{0}(\lambda)\chi_\infty(|x|) \\
&\Def  \Id + A_K + A_\infty.
\end{split}\]
As before,
\[
(P-\lambda)F(\Id - A_K - A_\infty + A_KA_\infty) = \Id -A_\infty A_K + A_\infty A_K A_\infty.
\]
Now
\[\|\chi F(\Id - A_K - A_\infty + A_KA_\infty)\chi\| \le C a(h),\]
and once again
\[\|A_\infty A_K\| = \Oh(h^\infty),\]
follows from propagation of singularities, this time in the following form (the proofs can again be found in, for example, \cite[Lemma 2]{d} and \cite[Lemma 5.1]{d-v1}): if $f \in L^2(\R^n)$ has compact support and $\|f\| \le h^{-N}$, then for all $\rho \in T^*\R^n$, 
\[
\rho \in \WFh\left(R_{0}(\lambda) f \right)\Longrightarrow \gamma_\rho^- \cap \WFh f \ne \varnothing.
\]
If further $\gamma_\rho^- \subset \{|x|>R_0\}$, then similarly
\[
\rho \in \WFh \left(R_W(\lambda) f \right)\Longrightarrow \gamma_\rho^- \cap \WFh f \ne \varnothing.
\]
\end{proof}

\end{document}